\documentclass[12pt,reqno,openbib,runningheads,a4paper]{amsart}%
\usepackage{eurosym}
\usepackage{amssymb}
\usepackage{amsfonts}
\usepackage{amsmath}
\usepackage{graphicx}
\usepackage[authoryear]{natbib}
\usepackage[legalpaper,bookmarks=true,colorlinks=true,linkcolor=blue,citecolor=blue]%
{hyperref}%
\providecommand{\U}[1]{\protect\rule{.1in}{.1in}}
\providecommand{\U}[1]{\protect\rule{.1in}{.1in}}
\newtheorem{theorem}{Theorem}[section]

\makeatletter
\renewcommand{\@biblabel}[1]{}
\makeatother
\begin{document}

\begin{center}
{\Large \textbf{Koml\'{o}s-Major-Tusn\'{a}dy approximations}}

{\Large \textbf{to increments of uniform empirical processes}}\medskip\medskip

{\large Abdelhakim Necir}$^{\ast}$\medskip

{\small \textit{Laboratory of Applied Mathematics, Mohamed Khider University,
Biskra, Algeria}}\medskip\medskip%
\[
\]

\end{center}

\noindent\textbf{Abstract}\medskip

\noindent The well-known Koml\'{o}s-Major-Tusn\'{a}dy inequalities [Z.
Wahrsch. Verw. Gebiete 32 (1975) 111-131; Z. Wahrsch. Verw. Gebiete 34 (1976)
33-58] provide sharp inequalities to partial sums of iid standard exponential
random variables by a sequence of standard Brownian motions. In this paper, we
employ these results to establish Gaussian approximations to weighted
increments of uniform empirical and quantile processes. This approach provides
rates to the approximations which, among others, have direct applications to
statistics of extreme values for randomly censored data.\medskip

\medskip


\vfill

\vfill

\noindent{\small $^{\text{*}}$Corresponding author:
\texttt{necirabdelhakim@yahoo.fr} \newline}

\section{\textbf{Introduction\label{sec1}}}

\noindent\cite{CsCsHM86} have constructed a probability space, denoted by
$\left(  \Omega,\mathcal{A},\mathbb{P}\right)  ,$ carrying a sequence of
independent random variables (rv's) $U_{1},U_{2},...$ uniformly distributed on
$\left(  0,1\right)  $ and a sequence of Brownian bridges $\left\{
B_{n}\left(  s\right)  ;0\leq s\leq1\right\}  _{n\geq1}$ such that for the
empirical process
\[
\alpha_{n}\left(  s\right)  :=\sqrt{n}\left(  G_{n}\left(  s\right)
-s\right)  ,\text{ }0\leq s\leq1
\]
and the quantile process%
\[
\beta_{n}\left(  s\right)  :=\sqrt{n}\left(  s-G_{n}^{-1}\left(  s\right)
\right)  ,\text{ }0\leq s\leq1
\]
where $G_{n}\left(  s\right)  :=n^{-1}\sum_{i=1}^{n}\boldsymbol{1}\left\{
U_{i}\leq s\right\}  $ and
\[
G_{n}^{-1}\left(  s\right)  :=\inf\left\{  t,\text{ }G_{n}\left(  s\right)
\geq s\right\}  ,0\leq s\leq1,
\]
with $G_{n}^{-1}\left(  0\right)  :=G_{n}^{-1}\left(  0+\right)  ,$ for
universal positive constants $a,$ $b$ and $c$%
\begin{equation}
\mathbb{P}\left\{  \sup_{0\leq s\leq d/n}\left\vert \beta_{n}\left(  s\right)
-B_{n}\left(  s\right)  \right\vert \geq n^{-1/2}\left(  a\log d+x\right)
\right\}  \leq be^{-cx},\label{ineq1}%
\end{equation}
for all  $0\leq x\leq d^{1/2}$ and $1\leq d\leq n,$ with the same inequality
holding for the supremum taken over $1-d/n\leq s\leq1.$ Thereby, they showed
that
\begin{equation}
\sup_{\lambda/n\leq s\leq1-\lambda/n}\frac{n^{\eta}\left\vert \beta_{n}\left(
s\right)  -B_{n}\left(  s\right)  \right\vert }{\left[  s\left(  1-s\right)
\right]  ^{1/2-\eta}}=O_{\mathbb{P}}\left(  1\right)  ,\label{approx1}%
\end{equation}
as $n\rightarrow\infty,$ for every fixed $0<\lambda<\infty$ and $0\leq
\eta<1/2,$ leading to%
\begin{equation}
\sup_{\lambda/n\leq s\leq1-\lambda/n}\frac{n^{\nu}\left\vert \alpha_{n}\left(
s\right)  -B_{n}\left(  s\right)  \right\vert }{\left[  s\left(  1-s\right)
\right]  ^{1/2-\nu}}=O_{\mathbb{P}}\left(  1\right)  ,\label{approx2}%
\end{equation}
for every fixed $0\leq\nu<1/4.$ The inequality $\left(  \ref{ineq1}\right)  $
is a result of Theorem 1.1 will  approximations $\left(  \ref{approx1}\right)
$ and $\left(  \ref{approx2}\right)  $ contain, respectively, in Theorem 2.1
and Corollary 2.1 of the above paper. Similar results may be found in
\cite{MZ-87}. These two Gaussian approximations remain powerful tools to
establish the asymptotic normality, among others, in statistics of extreme
values, see, e.g., \cite{cdm85} and \cite{peng2001}. In this paper, we are
concerned with Gaussian approximations of the increments%
\[
\alpha_{n}\left(  s;t\right)  :=\alpha_{n}\left(  t\right)  -\alpha_{n}\left(
t-s\right)  ,\text{ }0\leq s<t<1
\]
and%
\[
\beta_{n}\left(  s;t\right)  :=\beta_{n}\left(  t\right)  -\beta_{n}\left(
t-s\right)  ,\text{ }0\leq s<t<1.
\]
Such processes are used, for example, in goodness of fit test statistics
\citep[see,
e.g., Section 2 in][]{SW1982} and in nonparametric statistics for censored
data \citep[see,
e.g.,][]{DeEn96}. For convenience, we next use the notation $f\left(
s;t\right)  :=f\left(  t\right)  -f\left(  t-s\right)  ,$ $0\leq s<t<1,$ for
any measurable function $f.$ \cite{SW1982} (Theorem 1.2) showed that there
exist another Brownian bridge $\widetilde{B}\left(  s\right)  ,$ $0\leq
s\leq1$ such that
\begin{equation}
\sup_{cn^{-1}\log n\leq s<t}\frac{\left\vert \alpha_{n}\left(  s;t\right)
-\widetilde{B}\left(  s;t\right)  \right\vert }{s^{\nu}}=o_{\mathbb{P}}\left(
1\right)  ,\label{app}%
\end{equation}
for every $0<t<1$ and $c>0.$ By using Koml\'{o}s-Major-Tusn\'{a}dy
inequalities, \cite{CsCsHM86} (Theorem 4.6.1) also obtained a similar result
and proved that, in the probability space $\left(  \Omega,\mathcal{A}%
,\mathbb{P}\right)  ,$ we have%
\begin{equation}
\sup_{cn^{-1}\log n\leq s<t}\frac{\left\vert \alpha_{n}\left(  s;t\right)
-B_{n}\left(  s;t\right)  \right\vert }{s^{\nu}}=o_{\mathbb{P}}\left(
1\right)  ,\label{ap}%
\end{equation}
where $B_{n}$ is the same Brownian bridge as used in both approximations
$\left(  \ref{approx1}\right)  $ and $\left(  \ref{approx2}\right)  .$ The
authors are noticed, in their Remark 4.6.1, that $\left(  \ref{ap}\right)  $
is equivalent to%
\begin{equation}
\sup_{cn^{-1}\log n\leq s<t}\frac{\left\vert \beta_{n}\left(  s;t\right)
-B_{n}\left(  s;t\right)  \right\vert }{s^{\nu}}=o_{\mathbb{P}}\left(
1\right)  .\label{app2}%
\end{equation}
Otherwise, \cite{Alex87} (Remark 2.7, Assertion 2.7) gave a refinement of
$(\ref{app})$ and $(\ref{ap})$ to prove that there exists another Brownian
bridge $\widehat{B}\left(  s\right)  ,$ $0\leq s\leq1$ such that, for every
$0<\lambda<\infty$
\begin{equation}
\sup_{\lambda/n\leq s<t}\frac{\left\vert \alpha_{n}\left(  s;t\right)
-\widehat{B}\left(  s;t\right)  \right\vert }{s^{\nu}}=o_{\mathbb{P}}\left(
1\right)  .\label{app3}%
\end{equation}
It is worth mentioning, that the three Brownian bridges $B_{n},$
$\widetilde{B}$ and $\widehat{B}$ are not necessarily the same. Note also
that, for all large $n,$ $cn^{-1}\log n>n^{-1},$ then approximation $\left(
\ref{app3}\right)  $ is less restrictive and more useful than $\left(
\ref{ap}\right)  .$ But when we deal, for instance, to statistics of extreme
values for randomly censored data \citep[see,
e.g.,][]{BMN-2015} the rate of this approximation is needed. This, to our
knowledge, does not discussed yet in literature. In the following theorem we
answer to this issue by providing a new Gaussian approximation it term of a
sequence of Brownian bridges instead of their increments.

\begin{theorem}
\textbf{\label{Theo1}}On the probability space $\left(  \Omega,\mathcal{A}%
,\mathbb{P}\right)  ,$ carrying the sequence of iid rv's $U_{1},U_{2},...$
uniformly distributed on $\left(  0,1\right)  $ and the sequence of Brownian
bridges $B_{1},B_{2},...,$ for every $0<\lambda<\infty,$ $0\leq\eta<1/2,$ and
$0\leq\nu<1/4,$ we have approximations $\left(  \ref{approx1}\right)  $ and
$\left(  \ref{approx2}\right)  ,$ together with%
\begin{equation}
\sup_{\lambda/n\leq s<t}\frac{n^{\eta}\left\vert \beta_{n}\left(  s;t\right)
-B_{n}\left(  s\right)  \right\vert }{s^{1/2-\eta}}=O_{\mathbb{P}}\left(
1\right)  \label{approx3}%
\end{equation}
and
\begin{equation}
\sup_{\lambda/n\leq s<t}\frac{n^{\nu}\left\vert \alpha_{n}\left(  s;t\right)
-B_{n}\left(  s\right)  \right\vert }{s^{1/2-\nu}}=O_{\mathbb{P}}\left(
1\right)  .\bigskip\label{approx4}%
\end{equation}

\end{theorem}

\section{\textbf{Application to statistics for censored data\label{Sec}}}

\noindent Let $X_{1},...,X_{n}$ be $n\geq1$ independent copies of a
non-negative continuous random variable (rv) $X,$ defined over the probability
space with cumulative distribution function (cdf) $F.\ $These rv's are
censored to the right by a sequence of independent copies $Y_{1},...,Y_{n}$ of
a non-negative continuous rv $Y,$ independent of $X$ and having a cdf $G.$ At
each stage $1\leq j\leq n,$ we only can observe the rv's $Z_{j}:=\min\left(
X_{j},Y_{j}\right)  $ and $\delta_{j}:=\mathbf{1}\left\{  X_{j}\leq
Y_{j}\right\}  .$ If we denote by $H$ the cdf of the observed $Z^{\prime}s,$
then, in virtue of the independence of $X$ and $Y,$ we have $1-H=\left(
1-F\right)  \left(  1-G\right)  .$ We introduce two very crucial
sub-distribution functions $H^{\left(  i\right)  }\left(  z\right)
:=\mathbf{P}\left\{  Z_{1}\leq z,\delta_{1}=i\right\}  ,$ $i=0,1,$ for $z>0,$
so that one have $H\left(  z\right)  =H^{\left(  0\right)  }\left(  z\right)
+H^{\left(  1\right)  }\left(  z\right)  .$ The empirical counterparts are,
respectively, defined by%
\[
H_{n}^{\left(  0\right)  }\left(  z\right)  :=\frac{1}{n}\sum_{i=1}%
^{n}\mathbf{1}\left\{  Z_{i}\leq z\right\}  \left(  1-\delta_{i}\right)
,\text{ }H_{n}^{\left(  1\right)  }\left(  z\right)  :=\frac{1}{n}\sum
_{i=1}^{n}\mathbf{1}\left\{  Z_{i}\leq z\right\}  \delta_{i},
\]
and therefore
\[
H_{n}\left(  z\right)  :=n^{-1}\sum_{i=1}^{n}\mathbf{1}\left\{  Z_{i}\leq
z\right\}  =H_{n}^{\left(  0\right)  }\left(  z\right)  +H_{n}^{\left(
1\right)  }\left(  z\right)  .
\]
Let
\[
\xi_{i}:=\delta_{i}H^{\left(  1\right)  }\left(  Z_{i}\right)  +\left(
1-\delta_{i}\right)  \left(  \theta+H^{\left(  0\right)  }\left(
Z_{i}\right)  \right)  ,\text{ }i=1,...,n,
\]
be a sequence iid rv's uniformly distributed on $(0,1)$ \citep[][]{EnKo92},
and define the corresponding empirical cdf and empirical process by%
\[
\mathbb{U}_{n}\left(  s\right)  :=\frac{1}{n}\sum\limits_{i=1}^{n}%
\mathbf{1}\left\{  \xi_{i}\leq s\right\}  \text{ and }\alpha_{n}^{\ast}\left(
s\right)  :=\sqrt{n}\left(  \mathbb{U}_{n}\left(  s\right)  -s\right)  ,\text{
}0\leq s\leq1,
\]
respectively. Thereby we may represent, almost surely (a.s.), both
$H_{n}^{\left(  0\right)  }$ and $H_{n}^{\left(  1\right)  }$ in term of
$\mathbb{U}_{n},$ as follows $H_{n}^{\left(  0\right)  }\left(  v\right)
=\mathbb{U}_{n}\left(  H^{\left(  0\right)  }\left(  v\right)  +\theta\right)
-\mathbb{U}_{n}\left(  \theta\right)  ,$ for $0<H^{\left(  0\right)  }\left(
v\right)  <1-\theta,$ and $H_{n}^{\left(  1\right)  }\left(  v\right)
=\mathbb{U}_{n}\left(  H^{\left(  1\right)  }\left(  v\right)  \right)  ,$ for
$0<H^{\left(  1\right)  }\left(  v\right)  <\theta.$ For further details, see
for instance \cite{DeEn96}. From the previous representations, a.s., we may
write%
\[
\sqrt{n}\left(  \overline{H}_{n}^{\left(  0\right)  }\left(  v\right)
-\overline{H}^{\left(  0\right)  }\left(  v\right)  \right)  =-\alpha
_{n}^{\ast}\left(  1-\overline{H}^{\left(  0\right)  }\left(  v\right)
\right)  ,\text{ for }0<\overline{H}^{\left(  0\right)  }\left(  v\right)
<1-\theta
\]
and
\[
\sqrt{n}\left(  \overline{H}_{n}^{\left(  1\right)  }\left(  v\right)
-\overline{H}^{\left(  1\right)  }\left(  v\right)  \right)  =\alpha_{n}%
^{\ast}\left(  \overline{H}^{\left(  1\right)  }\left(  v\right)
;\theta\right)  ,\text{ for }0<\overline{H}^{\left(  1\right)  }\left(
v\right)  <\theta.
\]
By applying two approximations $\left(  \ref{approx2}\right)  $ and $\left(
\ref{approx4}\right)  ,$ there exists a sequence of Brownian bridges $\left\{
\mathcal{B}_{n}\left(  s\right)  ;\text{ }0\leq s\leq1\right\}  $ such that
for every $0<\lambda<\infty$ and $0\leq\xi<1/4,$%
\[
\sup_{\lambda/n\leq\overline{H}^{\left(  0\right)  }\left(  v\right)  \leq
1}\frac{n^{\xi}\left\vert \alpha_{n}^{\ast}\left(  1-\overline{H}^{\left(
0\right)  }\left(  v\right)  \right)  -\mathcal{B}_{n}\left(  1-\overline
{H}^{\left(  0\right)  }\left(  v\right)  \right)  \right\vert }{\left[
\overline{H}^{\left(  0\right)  }\left(  v\right)  \right]  ^{1/2-\xi}%
}=O_{\mathbb{P}}\left(  1\right)
\]
and%
\[
\sup_{\lambda/n\leq\overline{H}^{\left(  1\right)  }\left(  v\right)  <\theta
}\frac{n^{\xi}\left\vert \alpha_{n}^{\ast}\left(  \overline{H}^{\left(
1\right)  }\left(  v\right)  ;\theta\right)  -\mathcal{B}_{n}\left(
\overline{H}^{\left(  1\right)  }\left(  v\right)  \right)  \right\vert
}{\left[  \overline{H}^{\left(  1\right)  }\left(  v\right)  \right]
^{1/2-\xi}}=O_{\mathbb{P}}\left(  1\right)  .
\]
These approximations will be useful tools for asymptotic results to statistics
of extreme values for censored data, see for instance \cite{BMN-2015}.

\section{\textbf{Proof of Theorem }$\ref{Theo1}$\textbf{ \label{sec3}}}

\noindent\ Let $Y_{1}^{\left(  i\right)  },Y_{2}^{\left(  i\right)  },...,$
$\left(  i=1,2\right)  ,$ be two independent sequences of iid exponential rv's
with mean 1. From \cite{KMT75} inequalities, there exist two independent
copies $W^{\left(  i\right)  }\left(  z\right)  ,$ $0\leq z<\infty,$ $\left(
i=1,2\right)  ,$ of standard Brownian motion defined on a probability space,
such that for all real $x,$ we have
\begin{equation}
\mathbb{P}\left\{  \max_{1\leq k\leq m}\left\vert S_{k}^{\left(  i\right)
}-k-W^{\left(  i\right)  }\left(  k\right)  \right\vert \geq C\log
m+x\right\}  \leq Ke^{-\mu x},\label{KMT}%
\end{equation}
for $m=1,2,...,$ where $S_{k}^{\left(  i\right)  }:=\sum_{j=1}^{k}%
Y_{j}^{\left(  i\right)  },$ with $C,$ $K$ and $\mu$ are positive universal
constants independent of $i$ and $m.$ For each integer $n\geq2,$ we set
\[
Y_{j}\left(  n\right)  :=\left\{
\begin{array}
[c]{ll}%
Y_{\left[  n/2\right]  -j+1}^{\left(  1\right)  } & \text{for }j=1,...,\left[
n/2\right]  ,\medskip\\
Y_{n-j+2}^{\left(  2\right)  } & \text{for }j=\left[  n/2\right]  +1,...,n+1.
\end{array}
\right.
\]
Then $Y_{1}\left(  n\right)  ,Y_{2}\left(  n\right)  ,...,Y_{n+1}\left(
n\right)  $ are iid sequence of exponential rv's with mean 1. For further use,
we set $S_{m}\left(  n\right)  :=\sum\nolimits_{j=1}^{m}Y_{j}\left(  n\right)
,$ $m=1,...,n+1,$ and for the sake of notational simplicity, we will write
from now on, $S_{m}$ and $Y_{j}$ instead of $S_{m}\left(  n\right)  $ and
$Y_{j}\left(  n\right)  ,$ respectively, and will also use the usual
convention $S_{0}=0.$ It is easy to verify that, for each integer $n\geq2,$
the following process is a sequence of standard Brownian motions on $\left[
0,n+1\right]  :$%
\[
W_{n}\left(  s\right)  :=\left\{
\begin{array}
[c]{ll}%
W^{\left(  1\right)  }\left(  s\right)   & ,\text{ for }0\leq s\leq\left[
\dfrac{n}{2}\right]  ,\\
W^{\left(  1\right)  }\left(  \left[  \dfrac{n}{2}\right]  \right)
+W^{\left(  2\right)  }\left(  n+1-\left[  \dfrac{n}{2}\right]  \right)   & \\
\ \ \ \ \ \ \ \ \ \ \ \ \ \ \ -W^{\left(  2\right)  }\left(  n+1-s\right)   &
,\text{ for }\left[  \dfrac{n}{2}\right]  <s\leq n+1.
\end{array}
\right.
\]
Let us define the following two processes
\[
\widetilde{B}_{n}\left(  s\right)  :=n^{-1/2}\left(  sW_{n}\left(  n\right)
-W_{n}\left(  sn\right)  \right)  ,\text{ }0\leq s\leq1
\]
and
\[
\widetilde{\beta}_{n}\left(  s\right)  :=\sqrt{n}\left(  s-\widetilde
{U}_{\left[  sn\right]  :n}\right)  ,\text{ }0\leq s\leq1,
\]
where $\widetilde{U}_{k:n}:=S_{k}/S_{n+1},$ for $k=1,...,n,$ be a sequence of
the uniform order statistics, with the convention $\widetilde{U}_{0:n}%
=S_{0}\equiv0.$ We also define the uniform empirical process, corresponding to
$\widetilde{U}_{1:n},...,\widetilde{U}_{n:n},$ by%
\[
\widetilde{\alpha}_{n}\left(  s\right)  :=n^{1/2}\left(  \widetilde{G}%
_{n}\left(  s\right)  -s\right)  ,\text{ }0\leq s\leq1,
\]
where $\widetilde{G}_{n}\left(  s\right)  :=n^{-1}\sum_{i=1}^{n}%
\boldsymbol{1}\left\{  \widetilde{U}_{i:n}\leq s\right\}  .$ In their
inequalities $\left(  1.23\right)  $ and $\left(  1.24\right)  ,$
\cite{CsCsHM86}, showed that%
\[
\mathbb{P}\left\{  \sup_{0\leq s\leq d/n}\left\vert \widetilde{\beta}%
_{n}\left(  s\right)  -\widetilde{B}_{n}\left(  s\right)  \right\vert
\geq2n^{-1/2}\left(  a\log d+x\right)  \right\}
\]
and%
\[
\mathbb{P}\left\{  \sup_{1-d/n\leq s\leq1}\left\vert \widetilde{\beta}%
_{n}\left(  s\right)  -\widetilde{B}_{n}\left(  s\right)  \right\vert
\geq2n^{-1/2}\left(  a\log d+x\right)  \right\}  ,
\]
whenever $n_{0}<d<n$ and $0\leq x\leq d^{1/2}$ for suitably chosen positive
constants $n_{0},$ $a,$ $b$ and $c.$ Thereby they stated that%
\[
\sup_{\lambda/n\leq s\leq1-\lambda/n}\frac{n^{\eta}\left\vert \widetilde
{\beta}_{n}\left(  s\right)  -\widetilde{B}_{n}\left(  s\right)  \right\vert
}{\left[  s\left(  1-s\right)  \right]  ^{1/2-\eta}}=O_{\mathbb{P}}\left(
1\right)  =\sup_{\lambda/n\leq s\leq1-\lambda/n}\frac{n^{\nu}\left\vert
\widetilde{\alpha}_{n}\left(  s\right)  -\widetilde{B}_{n}\left(  s\right)
\right\vert }{\left[  s\left(  1-s\right)  \right]  ^{1/2-\nu}},
\]
for every $0<\lambda<\infty,$ $0\leq\eta<1/2,$ and $0\leq\nu<1/4.$ Next we
establish similar results to the increments
\[
\widetilde{\beta}_{n}\left(  s;t\right)  =\widetilde{\beta}_{n}\left(
t\right)  -\widetilde{\beta}_{n}\left(  t-s\right)  =\sqrt{n}\left(
s-\widetilde{U}_{\left[  nt\right]  :n}+\widetilde{U}_{\left[  n\left(
t-s\right)  \right]  :n}\right)  ,\text{ }0\leq s<t<1
\]
and%
\[
\widetilde{\alpha}_{n}\left(  s;t\right)  =\widetilde{\alpha}_{n}\left(
t\right)  -\widetilde{\alpha}_{n}\left(  t-s\right)  =\sqrt{n}\left(
\widetilde{G}_{n}\left(  s\right)  -\widetilde{G}_{n}\left(  t-s\right)
-s\right)  ,\text{ }0\leq s<t<1.
\]
To this end, we will follow similar steps as used for the proof of Theorem 1.1
(inequality 1.1) in \cite{CsCsHM86}. Let both $d$ and $n$ be sufficiently
large and $n_{0}>1,$ so that $n_{0}<d<n.$ For $a>0$ and $0\leq x\leq d^{1/2},$
we set%
\[
\mathbf{P}_{n}\left(  x;d\right)  :=\mathbb{P}\left\{  \sup_{0\leq s\leq
d/n}\left\vert \widetilde{\beta}_{n}\left(  s;t\right)  -\widetilde{B}%
_{n}\left(  s\right)  \right\vert \geq2n^{-1/2}\left(  a\log d+x\right)
\right\}  ,
\]
which is less than or equal to the sum of%
\[
\mathbf{P}_{1,n}\left(  x;d\right)  :=\mathbb{P}\left\{  \sup_{0\leq s\leq
d/n}\left\vert \widehat{\beta}_{n}\left(  s;t\right)  -\widetilde{\beta}%
_{n}\left(  s\right)  \right\vert \geq n^{-1/2}\left(  a\log d+x\right)
\right\}  ,
\]
and%
\[
\mathbf{P}_{2,n}\left(  x;d\right)  :=\mathbb{P}\left\{  \sup_{0\leq s\leq
d/n}\left\vert \widetilde{\beta}_{n}\left(  s\right)  -\widetilde{B}%
_{n}\left(  s\right)  \right\vert \geq n^{-1/2}\left(  a\log d+x\right)
\right\}  .
\]
Next we show that $\mathbf{P}_{1,n}\left(  x;d\right)  \overset{\mathbb{P}%
}{\rightarrow}0,$ as $n\rightarrow\infty.$ Indeed, let us write
\[
\left\vert \widehat{\beta}_{n}\left(  s;t\right)  -\widehat{\beta}_{n}\left(
s\right)  \right\vert =\sqrt{n}\left\vert \widetilde{U}_{\left[  nt\right]
:n}-\widetilde{U}_{\left[  n\left(  t-s\right)  \right]  :n}-\widetilde
{U}_{\left[  ns\right]  :n}\right\vert ,
\]
which equals%
\[
\sqrt{n}\frac{\left\vert S_{\left[  nt\right]  }-S_{\left[  n\left(
t-s\right)  \right]  }-S_{\left[  ns\right]  }\right\vert }{S_{n+1}}=\sqrt
{n}\frac{S_{\left\vert \left[  nt\right]  -\left[  n\left(  t-s\right)
\right]  -\left[  ns\right]  \right\vert }}{S_{n+1}},
\]
thus%
\[
\mathbf{P}_{1,n}\left(  x;d\right)  =\mathbb{P}\left\{  \frac{n}{S_{n+1}}%
\sup_{0\leq s\leq d/n}S_{\left\vert \left[  nt\right]  -\left[  n\left(
t-s\right)  \right]  -\left[  ns\right]  \right\vert }\geq a\log d+x\right\}
.
\]
Since $u\leq\left[  u\right]  \leq u+1,$ then it is easy to check that
\[
-2\leq\left[  nt\right]  -\left[  n\left(  t-s\right)  \right]  -\left[
ns\right]  \leq1,
\]
this implies that $\left\vert \left[  nt\right]  -\left[  n\left(  t-s\right)
\right]  -\left[  ns\right]  \right\vert \leq2,$ it follows that for $0\leq
s<t<1,$ we have $S_{\left\vert \left[  nt\right]  -\left[  n\left(
t-s\right)  \right]  -\left[  ns\right]  \right\vert }\leq S_{2},$ therefore
\[
\mathbf{P}_{1,n}\left(  x;d\right)  \leq\mathbb{P}\left\{  \frac{n}{S_{n+1}%
}S_{2}\geq a\log d+x\right\}  .
\]
By the law of large numbers, $\mathbb{P}\left\{  \left\vert n/S_{n+1}%
-1\right\vert \geq\epsilon\right\}  \rightarrow0,$ for any fixed
$0<\epsilon<1,$ this implies that
\[
\mathbf{P}_{1,n}\left(  x;d\right)  \leq\mathbb{P}\left\{  \left(
1-\epsilon\right)  S_{2}\geq a\log d+x\right\}  +\mathbb{P}\left\{  \left\vert
n/S_{n+1}-1\right\vert \geq\epsilon\right\}  .
\]
Note that $S_{2}$ is a sum of two iid standard exponential rv's, this means
that it follows the Gamma cdf with two parameters $\left(  2,1\right)  ,$ that
is $\mathbb{P}\left(  S_{2}>u\right)  =\left(  u+1\right)  e^{-u},$ therefore%
\[
\mathbb{P}\left\{  \left(  1-\epsilon\right)  S_{2}\geq a\log d+x\right\}
=\left(  \frac{a\log d+x}{1-\epsilon}+1\right)  \exp\left(  -\frac{a\log
d+x}{1-\epsilon}\right)  ,
\]
which tends to zero as $d\rightarrow\infty,$ hence $\mathbf{P}_{1,n}\left(
x;d\right)  \rightarrow0.$ On the other hand, from inequality $\left(
1.23\right)  $ in \cite{CsCsHM86}, we have $\mathbf{P}_{2,n}\left(
x;d\right)  \leq b\exp\left(  -cx\right)  ,$ thus $\mathbf{P}_{n}\left(
x;d\right)  \leq b\exp\left(  -cx\right)  ,$ too. Thereby, by using the latter
inequality with similar arguments as used for the proof of Theorem 2.1
(statement 2.2) of the same paper, we end up with
\begin{equation}
\sup_{\lambda/n\leq s<t}\frac{n^{\eta}\left\vert \widetilde{\beta}_{n}\left(
s;t\right)  -\widetilde{B}_{n}\left(  s\right)  \right\vert }{s^{1/2-\eta}%
}=O_{\mathbb{P}}\left(  1\right)  ,\label{beta-tild}%
\end{equation}
for every $0\leq\eta<1/2$ and $0<\lambda<\infty.$ Next we show that for every
$0<t<1$ and $0\leq\nu<1/4,$ we also have%
\begin{equation}
A_{n,\nu}\left(  t\right)  :=\sup_{\widetilde{U}_{1:n}\leq s<\widetilde
{U}_{t_{n}:n}}\frac{n^{\nu}\left\vert \widetilde{\alpha}_{n}\left(
s;t\right)  -\widetilde{B}_{n}\left(  s\right)  \right\vert }{s^{1/2-\nu}%
}=O_{\mathbb{P}}\left(  1\right)  ,\label{A}%
\end{equation}
where $t_{n}:=\left[  nt\right]  .$ Indeed, let us write%
\[
A_{n,\nu}\left(  t\right)  =\max_{1\leq k\leq t_{n}-1}\left\{  \sup
_{\widetilde{U}_{k:n}\leq s<\widetilde{U}_{k+1:n}}\frac{n^{\nu}\left\vert
\widetilde{\alpha}_{n}\left(  s;t\right)  -\widetilde{B}_{n}\left(  s\right)
\right\vert }{s^{1/2-\nu}}\right\}  ,
\]
and, for $0<\tau\leq1,$ set
\[
A_{n,\nu}\left(  t;\tau\right)  :=\max_{1\leq k\leq t_{n}-1}\left\{
\sup_{\widetilde{U}_{k:n}\leq s<\widetilde{U}_{k+1:n}}\frac{n^{\nu}\left\vert
\widetilde{\alpha}_{n}\left(  s;t\right)  -\widetilde{B}_{n}\left(  s\right)
\right\vert }{\left(  \tau k/n\right)  ^{1/2-\nu}}\right\}  .
\]
Observe that
\[
\left\{  \min_{1\leq k\leq n}\widetilde{U}_{k:n}/k\geq\tau\right\}
\subset\left\{  A_{n,\nu}\left(  t\right)  \leq\tau^{\nu-1/2}A_{n,\nu}\left(
t;1\right)  \right\}  ,
\]
and from assertion $\left(  2.9\right)  $ in \cite{CsCsHM86}, we have
\[
\mathbb{P}\left\{  \min_{1\leq k\leq n}\widetilde{U}_{k:n}/k\leq\tau\right\}
=\tau,
\]
it follows that
\[
\mathbb{P}\left\{  A_{n,\nu}\left(  t\right)  \leq\tau^{\nu-1/2}A_{n,\nu
}\left(  t;1\right)  \right\}  \geq1-\tau.
\]
Hence, to show that $A_{n,\nu}\left(  t\right)  =O_{\mathbb{P}}\left(
1\right)  ,$ it suffices to verify that $A_{n,\nu}\left(  t;1\right)
=O_{\mathbb{P}}\left(  1\right)  $ for sufficiently small $\tau.$ To this end,
we will first state that for $1\leq k\leq t_{n}-1$ and $\widetilde{U}%
_{k:n}\leq s<\widetilde{U}_{k+1:n},$ we have
\begin{equation}
\widetilde{\alpha}_{n}\left(  s;t\right)  -\widetilde{\beta}_{n}\left(
\frac{t_{n}}{n},\frac{k}{n}\right)  =O_{\mathbb{P}}\left(  n^{-1/2}\right)
.\label{equ}%
\end{equation}
Indeed, let us fix $\epsilon>0$ be small such that
\[
\widetilde{U}_{t_{n}+1:n}<\widetilde{U}_{t_{n}+2:n}-\epsilon<\widetilde
{U}_{t_{n}+2:n}+\epsilon<\widetilde{U}_{t_{n}+3:n},
\]
and set
\[
\mathcal{A}_{\epsilon,n}\left(  t\right)  :=\left\{  \left\vert \widetilde
{U}_{t_{n}+2:n}-t\right\vert <\epsilon\right\}  .
\]
Since $\widetilde{U}_{t_{n}+2:n}\overset{\mathbb{P}}{\rightarrow}t$ as
$n\rightarrow\infty$ then $\mathbb{P}\left(  \mathcal{A}_{\epsilon,n}\left(
t\right)  \right)  \downarrow1,$ as $n\rightarrow\infty.$ Hence, in the set
$\mathcal{A}_{\epsilon,n}\left(  t\right)  ,$ we have $\widetilde{U}%
_{t_{n}+2:n}-\epsilon<t<\widetilde{U}_{t_{n}+2:n}+\epsilon,$ which implies
that
\begin{equation}
\widetilde{U}_{t_{n}:n}<t<\widetilde{U}_{t_{n}+3:n}.\label{t}%
\end{equation}
Then, for $\widetilde{U}_{k:n}\leq s<\widetilde{U}_{k+1:n},$ we have
\[
\widetilde{G}_{n}\left(  t\right)  -\widetilde{G}_{n}\left(  t-s\right)
-s\geq\widetilde{G}_{n}\left(  \widetilde{U}_{t_{n}:n}\right)  -\widetilde
{G}_{n}\left(  \widetilde{U}_{t_{n}+3:n}-\widetilde{U}_{k:n}\right)
-\widetilde{U}_{k+1:n}.
\]
Note that
\[
\widetilde{U}_{t_{n}+3:n}-\widetilde{U}_{k:n}=\frac{S_{t_{n}+3}-S_{k}}%
{S_{n+1}}=\frac{S_{t_{n}-k+3}}{S_{n+1}}=\widetilde{U}_{t_{n}-k+3:n},
\]%
\[
\widetilde{G}_{n}\left(  \widetilde{U}_{t_{n}:n}\right)  =t_{n}/n\text{ and
}\widetilde{G}_{n}\left(  \widetilde{U}_{t_{n}-k+3:n}\right)  =\frac
{t_{n}-k+3}{n}.
\]
Then the right-side of the previous inequality is equal to%
\[
\frac{k-3}{n}-\widetilde{U}_{k+1:n}=\frac{k}{n}-\widetilde{U}_{k:n}-\left(
\widetilde{U}_{k+1:n}-\widetilde{U}_{k:n}\right)  -\frac{3}{n}.
\]
Note also $\widetilde{U}_{k+1:n}-\widetilde{U}_{k:n}=S_{1}/S_{n+1}$ and
$\widetilde{U}_{k:n}=\widetilde{U}_{t_{n}:n}-\widetilde{U}_{t_{n}-k:n},$ it
follows that%
\[
\widetilde{\alpha}_{n}\left(  s;t\right)  \geq\widetilde{\beta}_{n}\left(
\frac{k}{n},\frac{t_{n}}{n}\right)  -\sqrt{n}S_{1}/S_{n+1}-3/\sqrt{n}.
\]
By using the law of large numbers, we have with large probability
$n/S_{n+1}<2,$ then without loss of generality, we get%
\begin{equation}
\widetilde{\alpha}_{n}\left(  s;t\right)  \geq\widetilde{\beta}_{n}\left(
\frac{k}{n},\frac{t_{n}}{n}\right)  -3\left(  S_{1}+1\right)  /\sqrt
{n}.\label{i1}%
\end{equation}
Likewise, by using similar arguments as above we get
\[
\widetilde{G}_{n}\left(  t\right)  -\widetilde{G}_{n}\left(  t-s\right)
-s\leq\widetilde{G}_{n}\left(  \widetilde{U}_{t_{n}+3:n}\right)
-\widetilde{G}_{n}\left(  \widetilde{U}_{t_{n}:n}-\widetilde{U}_{k+1:n}%
\right)  -\widetilde{U}_{k:n},
\]
which implies that%
\begin{equation}
\widetilde{\alpha}_{n}\left(  s;t\right)  \leq\widetilde{\beta}_{n}\left(
\frac{k}{n},\frac{t_{n}}{n}\right)  +4/\sqrt{n}.\label{i2}%
\end{equation}
By letting $\zeta:=\max\left(  3\left(  S_{1}+1\right)  ,4\right)  ,$ the
inequalities $\left(  \ref{i1}\right)  $ and $\left(  \ref{i2}\right)  $
together give
\[
\left\vert \widetilde{\alpha}_{n}\left(  s;t\right)  -\widetilde{\beta}%
_{n}\left(  \dfrac{k}{n},\dfrac{t_{n}}{n}\right)  \right\vert <\zeta/\sqrt{n}.
\]
Since $\zeta=O_{\mathbb{P}}\left(  1\right)  ,$ hence $\widetilde{\alpha}%
_{n}\left(  s;t\right)  -\widetilde{\beta}_{n}\left(  \dfrac{k}{n}%
,\dfrac{t_{n}}{n}\right)  =O_{\mathbb{P}}\left(  n^{-1/2}\right)  $ which
meets  $\left(  \ref{equ}\right)  .$ It is clear that $A_{n,\nu}\left(
t;1\right)  $ is less than or equal to the sum of%
\[
L_{n}:=\max_{1\leq k\leq t_{n}-1}\left\{  \sup_{\widetilde{U}_{k:n}\leq
s<\widetilde{U}_{k+1:n}}\frac{n^{\nu}\left\vert \widetilde{\alpha}_{n}\left(
s;t\right)  -\widetilde{\beta}_{n}\left(  \dfrac{k}{n},\dfrac{t_{n}}%
{n}\right)  \right\vert }{\left(  k/n\right)  ^{1/2-\nu}}\right\}
\]
and%
\[
T_{n}:=\max_{1\leq k\leq t_{n}-1}\left\{  \sup_{\widetilde{U}_{k:n}\leq
s<\widetilde{U}_{k+1:n}}\frac{n^{\nu}\left\vert \widetilde{\beta}_{n}\left(
\dfrac{k}{n},\dfrac{t_{n}}{n}\right)  -\widetilde{B}_{n}\left(  s\right)
\right\vert }{\left(  k/n\right)  ^{1/2-\nu}}\right\}  .
\]
Making use of $\left(  \ref{equ}\right)  ,$ we infer that $L_{n}%
=O_{\mathbb{P}}\left(  1\right)  .$ Observe now that $T_{n}$ is less than or
equal to the sum of
\[
T_{n1}:=\max_{1\leq k\leq t_{n}-1}\left\{  \frac{n^{\nu}\left\vert
\widetilde{\beta}_{n}\left(  \dfrac{k}{n},\dfrac{t_{n}}{n}\right)
-\widetilde{B}_{n}\left(  \dfrac{k}{n}\right)  \right\vert }{\left(
k/n\right)  ^{1/2-\nu}}\right\}
\]
and%
\[
T_{n2}:=\max_{1\leq k\leq t_{n}-1}\left\{  \sup_{\widetilde{U}_{k:n}\leq
s<\widetilde{U}_{k+1:n}}\frac{n^{\nu}\left\vert \widetilde{B}_{n}\left(
s\right)  -\widetilde{B}_{n}\left(  \dfrac{k}{n}\right)  \right\vert }{\left(
k/n\right)  ^{1/2-\nu}}\right\}  .
\]
By letting $k/n=s$ and $t_{n}/n=t^{\ast},$ we may write%
\[
T_{n1}\leq\sup_{1/n\leq s<t^{\ast}}\frac{n^{\nu}\left\vert \widetilde{\beta
}_{n}\left(  s,t^{\ast}\right)  -\widetilde{B}_{n}\left(  s\right)
\right\vert }{s^{1/2-\nu}},
\]
which, by $\left(  \ref{beta-tild}\right)  ,$ is equal to $O_{\mathbb{P}%
}\left(  1\right)  .$ Let us now show that $T_{n2}=O_{\mathbb{P}}\left(
1\right)  $ too. To this end, we will follow similar procedures are those used
for the proof of assertion (2.21) in \cite{CsCsHM86}. Let us choose
$0<\nu<1/4$ and set $\delta:=\left(  1/4-\nu\right)  /2.$ For any $1\leq k\leq
n-1$ and $b\geq1,$ let $c_{k,n}^{\left(  \delta\right)  }:=k^{2\delta+1/2}/n$
and
\[
I_{k,n}\left(  b\right)  :=\left[  k/n-3bc_{k,n}^{\left(  \delta\right)
},k/n-3bc_{k,n}^{\left(  \delta\right)  }\right]
\]
and%
\[
D_{n,v}\left(  b\right)  :=\max_{1\leq k\leq t_{n}-1}\left\{  \sup_{s\in
I_{k,n}\left(  b\right)  }\frac{n^{\nu}\left\vert \widetilde{B}_{n}\left(
k/n\right)  -\widetilde{B}_{n}\left(  s\right)  \right\vert }{\left(
k/n\right)  ^{1/2-\nu}}\right\}  .
\]
Assertion $\left(  2.25\right)  $ in \cite{CsCsHM86} states that
\begin{equation}
\lim_{b\rightarrow\infty}\mathbb{P}\left\{  T_{n2}\geq D_{n,v}\left(
b\right)  \right\}  =0.\label{lim}%
\end{equation}
Then we have to show that $D_{n,v}\left(  b\right)  =O_{\mathbb{P}}\left(
1\right)  .$ Let us write%
\[
D_{n,v}\left(  b\right)  =\max_{1\leq k\leq t_{n}-1}\left\{  \sup
_{\mathbf{s}\in\left[  \mathbf{a}_{k}-h,\mathbf{a}_{k}+h\right]  \cap\left[
0,1\right]  }\frac{n^{\nu}\left\vert \widetilde{B}_{n}\left(  \mathbf{a}%
_{k}\right)  -\widetilde{B}_{n}\left(  \mathbf{s}\right)  \right\vert
}{\left(  k/n\right)  ^{1/2-\nu}}\right\}  ,
\]
where $\mathbf{a}_{k}=t-s,$ $\mathbf{s}=t^{\ast}-k/n$ and $h=\left\vert
\mathbf{a}_{k}\mathbf{-s}\right\vert .$ Since $\left\vert t-t^{\ast
}\right\vert \leq n^{-1}\leq bc_{k,n}^{\left(  \delta\right)  }$ and
$\left\vert s\mathbf{-}k/n\right\vert \leq3bc_{k,n}^{\left(  \delta\right)
},$ then $h\leq\left\vert t-t^{\ast}\right\vert +\left\vert s\mathbf{-}%
k/n\right\vert =4bc_{k,n}^{\left(  \delta\right)  }=:h_{\ast}.$ Let us write%
\[
D_{n,v}^{\ast}\left(  b\right)  :=\max_{1\leq k\leq t_{n}-1}\left\{
\sup_{\mathbf{s}\in\left[  \mathbf{a}_{k}-h_{\ast},\mathbf{a}_{k}+h_{\ast
}\right]  \cap\left[  0,1\right]  }\frac{n^{\nu}\left\vert \widetilde{B}%
_{n}\left(  \mathbf{a}_{k}\right)  -\widetilde{B}_{n}\left(  \mathbf{s}%
\right)  \right\vert }{\left(  k/n\right)  ^{1/2-\nu}}\right\}  .
\]
It is clear that $D_{n,v}^{\ast}\left(  b\right)  \leq D_{n,v}\left(
b\right)  .$ Hence, from $(\ref{lim}),$ it remain to show that $D_{n,v}^{\ast
}\left(  b\right)  =O_{\mathbb{P}}\left(  1\right)  .$ Indeed, for $d>0$
arbitrarily chosen, we have
\begin{align*}
&  \mathbb{P}\left\{  D_{n,v}^{\ast}\left(  b\right)  \geq d\left(  4b\right)
^{1/2}\right\}  \\
&  \leq\sum_{k=1}^{t_{n}-1}\mathbb{P}\left\{  \sup_{\mathbf{s}\in\left[
\mathbf{a}_{k}-h^{\ast},\mathbf{a}_{k}+h^{\ast}\right]  \cap\left[
0,1\right]  }\left\vert \widetilde{B}_{n}\left(  \mathbf{a}_{k}\right)
-\widetilde{B}_{n}\left(  \mathbf{s}_{k}\right)  \right\vert \geq d\left(
4b\right)  ^{1/2}k^{1/2-\nu}n^{-1/2}\right\}
\end{align*}
which may be rewritten into%
\begin{equation}
\sum_{k=1}^{t_{n}-1}\mathbb{P}\left\{  \sup_{\mathbf{s}\in\left[
\mathbf{a}_{k}-h^{\ast},\mathbf{a}_{k}+h^{\ast}\right]  \cap\left[
0,1\right]  }\left\vert \widetilde{B}_{n}\left(  \mathbf{a}_{k}\right)
-\widetilde{B}_{n}\left(  \mathbf{s}_{k}\right)  \right\vert \geq
dk^{1/4-\nu-\delta}h_{\ast}^{1/2}\right\}  .\label{sum}%
\end{equation}
From inequality (1.11) in \cite{CsCsHM86}, for a given Brownian bridge
$\mathbf{B}\left(  s\right)  ,$ $0\leq s\leq1$ defined on $\left(
\Omega,\mathcal{A},\mathbb{P}\right)  ,$ we have%
\[
\mathbb{P}\left\{  \sup_{\mathbf{s}\in\left[  \mathbf{a}-h,\mathbf{a}%
+h\right]  \cap\left[  0,1\right]  }\left\vert \mathbf{B}\left(
\mathbf{a}\right)  -\mathbf{B}\left(  \mathbf{s}\right)  \right\vert \geq
uh^{1/2}\right\}  \leq Au^{-1}\exp\left(  -u^{2}/8\right)  ,
\]
for any $0<\mathbf{a}<1,$ $h>0$ and $0<u<\infty,$ with a suitably chosen
universal constant $A.$ By applying this inequality we infer that $\left(
\ref{sum}\right)  $ is less than or equal to%
\[
\sum_{k=1}^{\infty}\frac{d^{-1}\exp\left(  -d^{2}k^{1/2-2\nu-2\delta
}/8\right)  }{k^{1/4-\nu-\delta}}=:\mathbf{P}\left(  d\right)  .
\]
Note that the series $\mathbf{P}\left(  d\right)  $ is uniformly convergent on
$d\geq1,$ and $\lim_{d\rightarrow\infty}\mathbf{P}\left(  d\right)  =0,$ it
follows that $D_{n,v}^{\ast}\left(  b\right)  =O_{\mathbb{P}}\left(  1\right)
,$ as sought. Let us now show that%
\[
\sup_{\lambda/n\leq s<t}\frac{n^{\nu}\left\vert \widetilde{\alpha}_{n}\left(
s;t\right)  -\widetilde{B}_{n}\left(  s\right)  \right\vert }{s^{1/2-\nu}%
}=O_{\mathbb{P}}\left(  1\right)  .
\]
Indeed, let us write%
\[
\sup_{\lambda/n\leq s<t}\frac{n^{\nu}\left\vert \widetilde{\alpha}_{n}\left(
s;t\right)  -\widetilde{B}_{n}\left(  s\right)  \right\vert }{s^{1/2-\nu}}\leq
A_{n,\nu}\left(  t\right)  +B_{n,\nu}\left(  t\right)  +C_{n,\nu}\left(
t\right)  +D_{n,\nu}\left(  t\right)  ,
\]
where $A_{n,\nu}\left(  t\right)  $ is that of $\left(  \ref{A}\right)  ,$
\[
B_{n,\nu}\left(  t\right)  :=\sup_{\left(  \lambda/n\right)  \wedge
\widetilde{U}_{1:n}\leq s<\widetilde{U}_{1:n}}\frac{n^{\nu}\left\vert
\widetilde{\alpha}_{n}\left(  s;t\right)  \right\vert }{s^{1/2-\nu}},
\]%
\[
C_{n,\nu}\left(  t\right)  :=\sup_{\widetilde{U}_{t_{n}:n}\leq s<\widetilde
{U}_{t_{n}:n}\vee t}\frac{n^{\nu}\left\vert \widetilde{\alpha}_{n}\left(
s;t\right)  -\widetilde{B}_{n}\left(  s\right)  \right\vert }{s^{1/2-\nu}}%
\]
and%
\[
D_{n,\nu}\left(  t\right)  :=\sup_{\left(  \lambda/n\right)  \wedge
\widetilde{U}_{1:n}\leq s<\widetilde{U}_{1:n}}\frac{n^{\nu}\left\vert
\widetilde{B}_{n}\left(  s\right)  \right\vert }{s^{1/2-\nu}}.
\]
It is clear that $B_{n,\nu}\left(  t\right)  $ is less than or equal to%
\[
\sup_{\left(  \lambda/n\right)  \wedge\widetilde{U}_{1:n}\leq s<\widetilde
{U}_{1:n}}\frac{n^{\nu+1/2}\left\vert G_{n}\left(  t\right)  -G_{n}\left(
t-s\right)  \right\vert }{s^{1/2-\nu}}+n^{\nu+1/2}\sup_{\left(  \lambda
/n\right)  \wedge\widetilde{U}_{1:n}\leq s<\widetilde{U}_{1:n}}s^{1/2+\nu}.
\]
Recall $\left(  \ref{t}\right)  $ and observe that for $\left(  \lambda
/n\right)  \wedge\widetilde{U}_{1:n}\leq s<\widetilde{U}_{1:n}$ (sufficiently
small), we have
\[
G_{n}\left(  t\right)  -G_{n}\left(  t-s\right)  \leq G_{n}\left(
\widetilde{U}_{t_{n}+3:n}\right)  -G_{n}\left(  \widetilde{U}_{t_{n}%
:n}-\widetilde{U}_{1:n}\right)  =\frac{t_{n}+3}{n}-\frac{t_{n}+1}{n}=\frac
{2}{n}%
\]
and%
\[
G_{n}\left(  t\right)  -G_{n}\left(  t-s\right)  \geq G_{n}\left(
\widetilde{U}_{t_{n}:n}\right)  -G_{n}\left(  \widetilde{U}_{t_{n}%
+3:n}\right)  =\frac{t_{n}}{n}-\frac{t_{n}+3}{n}=-\frac{3}{n}.
\]
It follows that $\left\vert G_{n}\left(  t\right)  -G_{n}\left(  t-s\right)
\right\vert \leq3/n,$ then it is easy to verify that%
\[
\sup_{\left(  \lambda/n\right)  \wedge\widetilde{U}_{1:n}\leq s<\widetilde
{U}_{1:n}}\frac{n^{\nu+1/2}\left\vert G_{n}\left(  t\right)  -G_{n}\left(
t-s\right)  \right\vert }{s^{1/2-\nu}}=O_{\mathbb{P}}\left(  1\right)  ,
\]
and since $n\widetilde{U}_{1:n}\overset{\mathbb{P}}{\rightarrow}1,$ we get%
\[
n^{\nu+1/2}\sup_{\left(  \lambda/n\right)  \wedge\widetilde{U}_{1:n}\leq
s<\widetilde{U}_{1:n}}s^{1/2+\nu}=\left(  n\widetilde{U}_{1:n}\right)
^{1/2+\nu}=O_{\mathbb{P}}\left(  1\right)  ,
\]
therefore $B_{n,\nu}\left(  t\right)  =o_{\mathbb{P}}\left(  1\right)  .$ For
the second term, we write%
\[
C_{n,\nu}\left(  t\right)  \leq\left(  \widetilde{U}_{t_{n}:n}\right)
^{1/2-\nu}\sup_{\widetilde{U}_{1:n}\leq s<\widetilde{U}_{n:n}}n^{\nu
}\left\vert \widetilde{\alpha}_{n}\left(  s;t\right)  -\widetilde{B}%
_{n}\left(  s\right)  \right\vert .
\]
In view of assertion $\left(  2.6\right)  $ of Theorem 2.2 in  \cite{CsCsHM86}%
, we have
\[
\sup_{0\leq s\leq1}\left\vert \widetilde{\alpha}_{n}\left(  s;t\right)
-\widetilde{B}_{n}\left(  s\right)  \right\vert =O\left(  \frac{\left(  \log
n\right)  ^{1/2}\left(  \log\log n\right)  ^{1/4}}{n^{1/4}}\right)  ,\text{
almost surely,}%
\]
it follows, since $0\leq\nu<1/4,$ that
\[
\sup_{\widetilde{U}_{1:n}\leq s<\widetilde{U}_{n:n}}n^{\nu}\left\vert
\widetilde{\alpha}_{n}\left(  s;t\right)  -\widetilde{B}_{n}\left(  s\right)
\right\vert =o_{\mathbb{P}}\left(  1\right)  .
\]
On the other hand $\widetilde{U}_{t_{n}:n}\overset{\mathbb{P}}{\rightarrow}t,$
it follows that $C_{n,\nu}\left(  t\right)  =o_{\mathbb{P}}\left(  1\right)
.$ We have $n\widetilde{U}_{1:n}\overset{\mathbb{P}}{\rightarrow}1,$ then it
is easy to show that $D_{n,\nu}\left(  t\right)  =O_{\mathbb{P}}\left(
1\right)  ,$ that we omit further details.  To summarize, we briefly stated
that
\begin{equation}
\sup_{\lambda/n\leq s\leq1-\lambda/n}\frac{n^{\eta}\left\vert \widetilde
{\beta}_{n}\left(  s\right)  -\widetilde{B}_{n}\left(  s\right)  \right\vert
}{\left[  s\left(  1-s\right)  \right]  ^{1/2-\eta}}=O_{\mathbb{P}}\left(
1\right)  =\sup_{\lambda/n\leq s\leq1-\lambda/n}\frac{n^{\nu}\left\vert
\widetilde{\alpha}_{n}\left(  s\right)  -\widetilde{B}_{n}\left(  s\right)
\right\vert }{\left[  s\left(  1-s\right)  \right]  ^{1/2-\nu}}.\label{ap1}%
\end{equation}
and showed that
\begin{equation}
\sup_{\lambda/n\leq s<t}\frac{n^{\eta}\left\vert \widetilde{\beta}_{n}\left(
s;t\right)  -\widetilde{B}_{n}\left(  s\right)  \right\vert }{s^{1/2-\eta}%
}=O_{\mathbb{P}}\left(  1\right)  =\sup_{\lambda/n\leq s<t}\frac{n^{\nu
}\left\vert \widetilde{\alpha}_{n}\left(  s;t\right)  -\widetilde{B}%
_{n}\left(  s\right)  \right\vert }{s^{1/2-\nu}}.\label{ap2}%
\end{equation}
On the other hand, we have for every $0<t<1,$%
\[
\left\{  \widetilde{\alpha}_{n}\left(  s\right)  ;0\leq s\leq1\right\}
\overset{\mathcal{D}}{=}\left\{  \alpha_{n}\left(  s\right)  ;0\leq
s\leq1\right\}
\]
and%
\[
\left\{  \widetilde{\beta}_{n}\left(  s\right)  ;0\leq s\leq1\right\}
\overset{\mathcal{D}}{=}\left\{  \beta_{n}\left(  s\right)  ;0\leq
s\leq1\right\}  .
\]
It follows that%
\[
\left\{  \widetilde{\alpha}_{n}\left(  s;t\right)  ;0\leq s<t\right\}
\overset{\mathcal{D}}{=}\left\{  \alpha_{n}\left(  s;t\right)  ;0\leq
s<t\right\}
\]
and
\[
\left\{  \widetilde{\beta}_{n}\left(  s;t\right)  ;0\leq s<t\right\}
\overset{\mathcal{D}}{=}\left\{  \beta_{n}\left(  s;t\right)  ;0\leq
s<t\right\}  ,
\]
with $\widetilde{B}_{n}$ is Brownian bridge for each $n.$ Then, having
established Gaussian approximations above one may construct a sequence
$U_{1},U_{2},...$ of iid rv's uniformly distributed on $\left[  0,1\right]  $
and a sequence of Brownian bridges $B_{1},B_{2},...$ defining on the
probability space $\left(  \Omega,\mathcal{A},\mathbb{P}\right)  $ such that
both $\left(  \ref{ap1}\right)  $ and $\left(  \ref{ap2}\right)  $ hold with
$\widetilde{\alpha}_{n},$ $\widetilde{\beta}_{n}$ and $\widetilde{B}_{n}$
replaced respectively by $\alpha_{n},$ $\beta_{n}$ and $B_{n}.$ This technique
for constructing a such probability space, described in Lemma 3.1.1 in
\cite{MCS83}, is used for instance in both \cite{CsCsHM86} and \cite{MZ-87}.

\end{document}